\documentclass[11pt,reqno]{amsart}
\usepackage{amssymb,amsmath,amsthm,amsfonts}
\usepackage[english]{babel}
\usepackage{mathrsfs,a4wide,dsfont}
\usepackage[utf8]{inputenc}
\usepackage{color}
\usepackage{comment}
\usepackage{bm}

\theoremstyle{plain}
\newtheorem{theorem}{Theorem}[section]

\newtheorem{corollary}[theorem]{Corollary}
\newtheorem{remark}[theorem]{Remark}

\theoremstyle{definition}
\theoremstyle{remark}
\numberwithin{equation}{section}

\newcommand{\R}{{\mathbb R}}

\newcommand{\N}{{\mathbb N}}

\newcommand{\e}{\varepsilon}

\definecolor{dg}{rgb}{0.01, 0.75, 0.24}

\title[Asymptotic analysis of a family of non-local functionals on sets]{Asymptotic analysis of a family of non-local functionals on sets}

\author[M.\,Eleuteri]
{Michela Eleuteri}
\address[M.\,Eleuteri]{Dipartimento di Scienze Fisiche, Informatiche e Matematiche, Universit\`a degli Studi di Modena e Reggio Emilia, via Campi 213, 41125 Modena, Italy.}
\email{michela.eleuteri@unimore.it}

\author[L.\,Lussardi]
{Luca Lussardi}
\address[L.\,Lussardi]{Dipartimento di Scienze Matematiche ``G.L.\,Lagrange'', Politecnico di Torino, c.so Duca degli Abruzzi 24, 10129 Torino, Italy.}
\email[]{luca.lussardi@polito.it}

\author[A.\,Torricelli]
{Andrea Torricelli}
\address[A.\,Torricelli]{Dipartimento di Scienze Fisiche, Informatiche e Matematiche, Universit\`a degli Studi di Modena e Reggio Emilia, via Campi 213, 41125 Modena, Italy.}
\email{andrea.torricelli@unimore.it}

\begin{document}

\baselineskip3.4ex

\vspace{0.5cm}
\begin{abstract}
{\small We study the asymptotic behavior of a family of functionals which penalize a short-range interaction of convolution type between a finite perimeter set and its complement. We first compute the pointwise limit and we obtain a lower estimate on more regulars sets. Finally, some examples are discussed. 
}
\end{abstract}

\maketitle

\section{Introduction}

In this paper we study the asymptotic behavior, as $\e \to 0$, of the family of functionals 
\[
\mathscr F_\e(E)=\frac{1}{\e}\int_{E^c \cap \Omega}f(G_\e \ast \chi_{E\cap \Omega})\,dx.
\]
where $\Omega\subset \R^N$ is open and bounded, $N>1$, $E$ is a set of finite perimeter in $\Omega$, $f$ is given and $G_\e(z)=\frac{1}{\e^N}G(\frac{z}{\e})$ where $G$ is a suitable kernel. Our analysis has been inspired by a paper by Miranda et al.\,in \cite{MPPP} where the case $f(t)=t$ is considered and $G$ is the Gauss-Weierstrass kernel, namely $G(z)=\frac{1}{(4\pi)^{N/2}}e^{-|z|^2/4}$ (see also \cite{GB} and \cite{L} for smoother sets and \cite{AB} for similar convolution approximation). More precisely, in \cite{MPPP} it is proven that the pointwise limit is, up to a constant, the perimeter of $E$ in $\Omega$. A more general kernel $G$ has been investigated, in the context of optimal partition problems, by Esedo\=glu and Otto \cite{EO} where $G$ is smooth and non-negative, radially symmetric and satisfying the following conditions:
\[
\int_{\R^N} G(z)\,dz=1, \quad \int_{\R^N}|z|G(z)\,dz<+\infty, \quad |\nabla G(z)|\lesssim G\left(\frac{z}{2}\right), \quad \nabla G(z)\cdot z \le 0.
\]
On the other hand, as in \cite{MPPP}, Esedo\=glu and Otto consider only the case $f(t)=t$, but they prove a complete $\Gamma$-convergence result for the family $\{\mathscr F_\e\}_{\e>0}$ on finite perimeter sets with respect to the strong $L^1$-convergence. A very similar result has been obtained more recently by Berendsen and Pagliari \cite{BP}. As far as we know, the last result is due to Pagliari \cite{P} where he essentially remove the radial symmetry of $G$ and he obtain, as limit, an anisotropic perimeter. 

In this paper we try to investigate the general situation. We assume that $G$ is even, non-negative, supported on the unit closed ball and with $\int_{\R^N} G(z)\,dz=1$. First of all we are able to compute the pointwise limit, as $\e\to0$, of $\mathscr F_\e(E)$ whenever $f$ is $C^1$, non-decreasing and $f(0)=0$. It turns out (see Theorem \ref{main1}) that for any $E \subset \R^N$ with finite perimeter in $\Omega$
\[
\lim_{\e\to0}\mathscr F_\e(E)=\int_{\partial^*E \cap \Omega}\int_0^1f\left(\int_{\{z \cdot \nu_E(x)\ge t\}}G(z)\,dz\right)\,dt\,d\mathcal H^{N-1}(x)
\]
where $\partial^*E$ is the reduced boundary of $E$ and $\nu_E(x)$ is the outer unit normal at $E$. In view to have a $\Gamma$-convergence result we investigate also the lower estimate. Unfortunately, the technique of Esedo\=glu and Otto \cite{EO} does not work in our situation: it is crucial for them to switch the order of integration, that is impossible for us since we have $f$ between the exterior integral and the convolution one. It seems that this difficulty cannot be easily overcome in the general situation. We are able to show (see Theorem \ref{main2}) a $\Gamma$-liminf inequality only on graphs of $C^1$ functions with respect to the $C^1$-uniform convergence. Actually, it is easy to generalize such a inequality in the case of sets which are locally graphs of  $C^1$ functions with respect to a suitable convergence (see Remark \ref{rem1}). Finally, we also prove (see Theorem \ref{main3}) that if $f$ is also convex, then the pointwise limit is lower semicontinuous with respect to the strong $L^1$-convergence, which suggests that for $f$ convex the $\Gamma$-limit in the strong $L^1$-convergence should be the pointwise limit. At the end of the paper we will also discuss some examples.

 \section{Notation and preliminaries}

\subsection{Notation}

In what follows $N\in \N$ with $N\ge 1$. For any $r>0$ and $x\in \R^N$ the notation $B^d_r(x)$ stands for the open ball in $\R^d$ centered at $x$ with radius $r$, while $\mathbb S^{N-1}=\partial B^N_1(0)$. If $A\subseteq \R^N$ we also denote by $\mathcal H^k(A)$ the Hausdorff measure of $A$ of dimension $k\in \{0,1,\dots,N\}$ ($\mathcal H^0$ is the counting measure). If $A_h,A$ are measurable subsets of $\R^N$,  then $A_h \to A$ in $L^1(\R^N)$ (or $L^1_{\rm loc}(\R^N)$) means that $\chi_{A_h}\to \chi_A$ in $L^1(\R^N)$ (respectively $L^1_{\rm loc}(\R^N)$). Finally, for any $A\subseteq \R^N$ we let $A^c=\R^N\setminus A$.

\subsection{Finite perimeter sets}

We recall some notion on finite perimeter sets in euclidean space; for details we refer to \cite{AFP-libro}. Let $\Omega$ be an open  subset of $\R^N$. A measurable set $E \subseteq \R^N$ is said to be a {\it set of finite perimeter in $\Omega$} if 
\[
\mathcal P(E,\Omega)=\sup\left\{\int_E{\rm div}\,X(x)\,dx: X \in C_c^1(\Omega;\R^N),\,\|X\|_{\infty}\leq 1 \right\}<+\infty.
\]
The quantity $\mathcal P(E,\Omega)$ is called {\it perimeter of $E$ in $\Omega$}. Finite perimeter sets have nice boundary in a measure theoretical sense. Precisely, one can define a subset of $E$ as the set of points $x$ where there exists a unit vector $\nu_E(x)$ such that:
\begin{equation}\label{conv}
\frac{x-E}{r} \to \{y \in \R^N: y\cdot \nu_E(x) \ge 0\},\,\text{ in $L^1_{\rm loc}(\R^N)$ as $r \to 0$},
\end{equation}
and which is referred to as the {\it outer normal to $E$ at $x$}. The set where $\nu_E(x)$ exists is called the \textit{reduced boundary of $E$} and is denoted by $\partial^*E$. It turns out that, for any $E$ set of finite perimeter in $\Omega$, we have $\mathcal P(E,\Omega)=\mathcal H^{N-1}(\partial^*E \cap \Omega)$. The reduced boundary of $E$ plays the role of the topological boundary also in the sense of the integration by parts. Indeed, one can show that, if $E$ is a set of finite perimeter in $\Omega$, then the following Gauss-Green formula holds true:
\begin{equation}\label{GG}
\int_E {\rm div}\,X(x)\,dx=\int_{\partial^*E}X(x) \cdot \nu_E(x)\,d\mathcal H^{N-1}(x), \quad \forall X \in C_c^1(\Omega;\R^N). 
\end{equation}
Finite perimeter sets satisfy good properties for the Calculus of Variations: for instance, if $E_h,E$ have finite perimeter in $\Omega$ and $E_h \stackrel{L^1}{\to}E$, then 
\[
\mathcal P(E,\Omega) \le \liminf_{h\to+\infty}\mathcal P(E_h,\Omega).
\]

\section{Setting of the problem and main results}

Let $N>1$, let $G \colon \R^N \to [0,+\infty)$ be of class $C^\infty$ such that
\[
{\rm supp}\,G=\overline{B^N_1(0)}, \quad G(-x)=G(x),\quad  \int_{\R^N} G(x)\,dx=1. 
\]
For any $\e>0$ and for any $x\in \R^N$, let 
\[
G_{\varepsilon}(x)=\frac{1}{\varepsilon^N}G\left(\frac{x}{\varepsilon}\right).
\]
We consider a continuous and non-decreasing function $f\colon [0,+\infty) \to\R$ with $f(0)=0$. Let $\Omega \subset \R^N$ be open bounded. We denote by $\mathcal P_N(\Omega)$ the set of all sets of finite perimeter in $\Omega$.  For any $\e>0$, we introduce the functional $\mathscr F_\e \colon \mathcal P_N(\Omega) \to [0,+\infty)$ defined by
\begin{equation}
\label{Fepsilon}
\mathscr F_\e(E)=\frac{1}{\e}\int_{E^c \cap \Omega}f(G_\e \ast \chi_{E\cap \Omega})\,dx.
\end{equation}
In order to state our main results, we introduce the function $\theta \colon \mathbb S^{N-1}\to [0,+\infty)$ given by
\begin{equation}\label{theta}
\theta(\nu)=\int_0^1 f\left(\int_{\{x \cdot \nu\ge t\}}G(x)\,dx\right)\,dt.
\end{equation}
Let $\mathscr F \colon \mathcal P_N(\Omega) \to [0,+\infty)$ be the functional given by 
\[
\mathscr F(E)=\int_{\partial^*E \cap \Omega}\theta(\nu_E(x))\,d\mathcal H^{N-1}(x).
\]
Our first main result concerns the pointwise limit of $\mathscr F_\e$ on $\mathcal P_N(\Omega)$. 

\begin{theorem}{\bf(Pointwise limit)}\label{main1}
Assume $f$ of class $C^1$. Let $E\in \mathcal P_N(\Omega)$. Then
\[
\lim_{\e\to0}\mathscr F_\e(E)=\mathscr F(E).
\]
\end{theorem}

On the other hand, we are also able to prove a lower estimate on graphs.

\begin{theorem}{\bf(Lower estimate)}\label{main2}
Let $D \subset \mathbb{R}^{N-1}$ be open and bounded with Lipschitz boundary, let $u_h,u \in C^{1,1}(D)$, with $u_h, u > 0$ on $D$ such that $u_h\to u$ uniformly in $C^1(D)$. Let $E_h,E$ be given by
\[
E_h = \{(x,y) \in \mathbb{R}^{N-1} \times \R: x \in D,\,0 \le y \le u_h(x)\},
\]
\[
E = \{(x,y) \in \mathbb{R}^{N-1} \times \R: x \in D,\,0 \le y \le u(x)\}.
\]
Then, for any positive infinitesimal sequence $(\e_h)$ it holds
\[
\liminf_{h \rightarrow +\infty} \mathscr F_{\e_h}(E_h) \ge \, \mathscr F(E).
\]
\end{theorem}

\begin{remark}\label{rem1}
It is not difficult to see that Theorem \ref{main2} can be generalized to uniformly $C^{1,1}$-regular sets in $\Omega$ with respect to a suitable notion of uniform convergence. Precisely, a set $E\subset \R^N$ is said to be {\rm uniformly $C^{1,1}$-regular set in $\Omega$} if there exist $L,\delta>0$ such that for every $x \in \partial E \cap \Omega$ there exist $D^x \subseteq \R^{N-1}$ open and a function $u^x \in C^{1,1}(D^x)$ such that:
\begin{itemize}
\item $\partial E \cap \Omega \cap B^N_\delta(x)$ is the graph of $u^x$;
\item $\|\nabla u^x\|_{\infty}\leq L$.
\end{itemize}
On the set of all uniformly $C^{1,1}$-regular sets in $\Omega$ we put a convergence of sequences. Precisely, we say that $E_h$ converges to $E$ if there exist $\delta,L>0$ such that for every $x \in \partial E \cap \Omega$ there exist $D^x \subseteq \R^{N-1}$ open and functions $u_h^x,u^x \in C^{1,1}(D^x)$ such that:
\begin{itemize}
\item $\partial E_h \cap \Omega \cap B^N_\delta(x),\partial E \cap \Omega \cap B^N_\delta(x)$ are the graphs of $u_h^x,u^x$ respectively; 
\item $\|\nabla u_h^x\|_{\infty}\leq L$ and $\|\nabla u^x\|_{\infty}\leq L$;
\item $u_h^x \to u^x$ uniformly in $C^1(D^x)$.
\end{itemize}
It is easy to see that with respect to this type of convergence the lower estimate 
\[
\liminf_{h \rightarrow +\infty} \mathscr F_{\e_h}(E_h) \ge \mathscr F(E)
\]
follows as a simple consequence of Theorem \ref{main2}.
\end{remark}

Combining Theorem \ref{main1} with Theorem \ref{main2} and Remark \ref{rem1} we eventually obtain a $\Gamma$-convergence result.

\begin{corollary}
The family $\{\mathcal F_\e\}_{\e>0}$ $\Gamma$-converges to $\mathcal F$ as $\e\to 0$ on uniformly $C^{1,1}$-regular sets with respect to the convergence introduced in Remark \ref{rem1}. 
\end{corollary}

\begin{remark}
We do not expect compactness of equibounded sequences of uniformly $C^{1,1}$-regular sets. Nevertheless, at least if $f(t)\ge mt$ for some $m>0$, equibounded sequences are compact in $L^1$. Indeed, if $(\e_h)$ is a positive and infinitesimal sequence and $(E_h)$ be a sequence in $\mathcal P_N(\Omega)$ with $\mathscr F_{\e_h}(E_h)\le c$ for some $c\ge 0$, we get 
\[
c\ge \mathscr F_{\e_h}(E_h)\ge \frac{m}{\e}\int_{E^c \cap \Omega}G_\e \ast \chi_{E\cap \Omega}\,dx
\]
and the compactness follows by \cite[Lemma A.4]{EO} (see also \cite[Thm.\,3.1]{AB}).
\end{remark}

The next and last result suggests that the $\Gamma$-limit on $\mathcal P_N(\Omega)$ of the family $(\mathscr F_\e)_{\e>0}$ with respect to the $L^1$-convergence could be really $\mathscr F$, at least if $f$ is convex.

\begin{theorem}\label{main3}
If $f$ is convex then the functional $\mathscr F \colon \mathcal P_N(\Omega) \to \R$ is lower semicontinuous with respect to the $L^1$-topology. 
\end{theorem}

\section{The pointwise limit}

In this section we prove Theorem \ref{main1}. The main idea comes from the technique used in \cite{MPPP}. We divide the proof in some steps.
\\
\\
{\sc Step 1}. We claim that for any $E\in \mathcal P_N$ we have 
\begin{equation}\label{claim1}
\mathscr F_\e(E)=\frac{1}{\e}\int_{\partial^*E}\int_0^\e X(\eta,x)\cdot \nu_E(x)\,d\eta\,d\mathcal H^{N-1}(x),
\end{equation}
where for any $\eta>0$ and for any $x\in \partial^*E$ 
\[
X(\eta,x)=\frac{1}{\eta^N}\int_{E^c}f'(G_\eta \ast \chi_E(y))G\left(\frac{y-x}{\eta}\right)\frac{y-x}{\eta}\,dy.
\]
For any $\eta>0$ and any $y\in \R^N$ we have, using the Gauss-Green formula \eqref{GG},  
\[
\begin{aligned}
\frac{d}{d\eta}&f(G_\eta\ast \chi_E(y))\\
&=-f'(G_\eta \ast \chi_E(y))\frac{1}{\eta^{N+1}}\int_{\R^N}\left(NG\left(\frac{y-x}{\eta}\right)+\nabla G\left(\frac{y-x}{\eta}\right)\cdot \frac{y-x}{\eta}\right)\chi_E(x)\,dx\\
&=f'(G_\eta \ast \chi_E(y))\frac{1}{\eta^N}\int_{\R^N}{\rm div}_x\left(G\left(\frac{y-x}{\eta}\right)\frac{y-x}{\eta}\right)\chi_E(x)\,dx\\
&=f'(G_\eta \ast \chi_E(y))\frac{1}{\eta^N}\int_{\partial^*E}G\left(\frac{y-x}{\eta}\right)\frac{y-x}{\eta}\cdot \nu_E(x)\,d\mathcal H^{N-1}(x).
\end{aligned}
\]
Now notice that, since $G_\e\ast \chi_E\to \chi_E$ in $L^1(\R^N)$ as $\e\to0$, we can say that for any $y\in \R^N$ 
\[
f(G_\e\ast \chi_E(y))-f(\chi_E(y))=\int_0^\e\frac{d}{d\eta}f(G_\eta\ast \chi_E(y))\,d\eta,
\]
from which we get, using the fact that $f(0)=0$, 
\[
\begin{aligned}
\mathscr F_\e(E)&=\frac{1}{\e}\int_{E^c}f(G_\e\ast \chi_E(y))-f(\chi_E(y))\,dy\\
&=\frac{1}{\e}\int_{E^c} \int_0^\e \frac{d}{d\eta}f(G_\eta\ast \chi_E(y))\,d\eta\,dy\\
&=\frac{1}{\e}\int_{\partial^*E}\int_0^\e\frac{1}{\eta^N} \int_{E^c}f'(G_\eta \ast \chi_E(y))G\left(\frac{y-x}{\eta}\right)\frac{y-x}{\eta}\,dy\,d\eta \cdot \nu_E(x)\,d\mathcal H^{N-1}(x)\\
&=\frac{1}{\e}\int_{\partial^*E}\int_0^\e X(\eta,x)\cdot \nu_E(x)\,d\eta\,d\mathcal H^{N-1}(x)
\end{aligned}
\]
hence \eqref{claim1}. 
\\
\\
{\sc Step 2}. We claim that for any $x\in \partial^*E$ we have 
\begin{equation}\label{claim2}
\lim_{\e\to0}X(\e,x)=\int_{\{z\cdot \nu_E(x)\ge 0\}}f'\left(\int_{\{(v-z)\cdot \nu_E(x)\ge 0\}}G(v)\,dv\right)G(z)z\,dz.
\end{equation}
First of all we have 
\[
\begin{aligned}
X(\e,x)&=\frac{1}{\e^N}\int_{E^c}f'(G_\e \ast \chi_E(y))G\left(\frac{y-x}{\e}\right)\frac{y-x}{\e}\,dy\\
&=\frac{1}{\e^N}\int_{E^c}f'\left(\frac{1}{\e^N}\int_EG\left(\frac{y-w}{\e}\right)\,dw\right)G\left(\frac{y-x}{\e}\right)\frac{y-x}{\e}\,dy.
\end{aligned}
\]
Performing first the change of variable $y=x+\e z$ and then $w=x+\e z-\e v$, we obtain
\[
\begin{aligned}
X(\e,x)&=\int_{\frac{E^c-x}{\e}}f'\left(\frac{1}{\e^N}\int_EG\left(\frac{x+\e z-w}{\e}\right)\,dw\right)G(z)z\,dz\\
&=\int_{\frac{E^c-x}{\e}}f'\left(\int_{\frac{x-E}{\e}+z}G(v)\,dv\right)G(z)z\,dz.
\end{aligned}
\]
Passing to the limit as $\e\to 0$ using \eqref{conv} and applying the Dominated convergence Theorem we easily get \eqref{claim2}. 
\\
\\
{\sc Step 3}. We claim that for any $x\in \partial^*E$ it holds 
\begin{equation}\label{claim3}
\int_{\{z\cdot \nu_E(x)\ge 0\}}f'\left(\int_{\{(v-z)\cdot \nu_E(x)\ge 0\}}G(v)\,dv\right)G(z)z\,dz\cdot \nu_E(x)=\theta(\nu_E(x)).
\end{equation}
First of all observe any $z \in \R^N$ with $z\cdot \nu_E(x)\ge 0$ can be written in a unique way as $z=\bar z+ t\nu_E(x)$ with $\bar z \cdot \nu_E(x)=0$ and $t\geq 0$. In particular, $z\cdot \nu_E(x)=(\bar z+t\nu_E(x))\cdot \nu_E(x)=t$. Moreover, since $G$ is supported on $\overline{B_1^N(0)}$, we can consider $t\in [0,1]$ obtaining 
\[
\begin{aligned}
\int_{\{z\cdot \nu_E(x)\ge 0\}}&f'\left(\int_{\{(v-z)\cdot \nu_E(x)\ge 0\}}G(v)\,dv\right)G(z)z\,dz \cdot \nu_E(x)\\
&=\int_{\{z\cdot \nu_E(x)\ge 0\}}f'\left(\int_{\{(v-z)\cdot \nu_E(x)\ge 0\}}G(v)\,dv\right)G(z)z\cdot \nu_E(x)\,dz\\
&=\int_0^1\int_{\{\bar z\cdot \nu_E(x)=0\}}f'\left(\int_{\{v\cdot \nu_E(x)\ge t\}}G(v)\,dv\right)G(\bar z+t\nu_E(x))\,t\,dt\,d\bar z\\
&=\int_0^1f'\left(\int_{\{v\cdot \nu_E(x)\ge t\}}G(v)\,dv\right)\int_{\{\bar z\cdot \nu_E(x)=0\}}G(\bar z+t\nu_E(x))\,t\,d\bar z\,dt\\
&=\int_0^1f'\left(\int_{\{v\cdot \nu_E(x)\ge t\}}G(v)\,dv\right)\int_{\{z\cdot \nu_E(x)=t\}}G(z)\,d\mathcal H^{N-1}(z)\,t\,dt.
\end{aligned}
\]
Finally, we remark that  
\[
\begin{aligned}
\frac{d}{dt}&\int_{\{v\cdot \nu_E(x)\ge t\}}G(v)\,dv\\
&=\lim_{h\to0}\frac{1}{h}\left(\int_{\{v\cdot \nu_E(x)\ge t+h\}}G(v)\,dv-\int_{\{v\cdot \nu_E(x)\ge t\}}G(v)\,dv\right)\\
&=-\lim_{h\to0}\frac{1}{h}\int_{\{t\le v \cdot \nu_E(x) \le t+h\}}G(v)\,dv\\
&=-\int_{\{v\cdot \nu_E(x)=t\}}G(v)\,dv.
\end{aligned}
\]
Integrating by parts we finally get 
\[
\begin{aligned}
\int_0^1&f'\left(\int_{\{v\cdot \nu_E(x)\ge t\}}G(v)\,dv\right)\int_{\{z\cdot \nu_E(x)=t\}}G(z)\,d\mathcal H^{N-1}(z)\,t\,dt\\
&=-\int_0^1\frac{d}{dt}f\left(\int_{\{v\cdot \nu_E(x)\ge t\}}G(v)\,dv\right)\,t\,dt\\
&=-f\left(\int_{\{v\cdot \nu_E(x)\ge t\}}G(v)\,dv\right)\,t\bigg|_0^1+\int_0^1f\left(\int_{\{v\cdot \nu_E(x)\ge t\}}G(v)\,dv\right)\,dt\\
&=\theta(\nu_E(x))
\end{aligned}
\]
where $\theta$ has been introduced in \eqref{theta}. This concludes the proof of \eqref{claim3}.
\\
\\
{\sc Step 4}. We easily conclude. Using \eqref{claim1}, \eqref{claim2}, \eqref{claim3}, De l'H\^opital rule and the Dominated convergence Theorem we deduce that
\[
\begin{aligned}
\lim_{\e\to0}\mathscr F_\e(E)&=\lim_{\e\to0}\frac{1}{\e}\int_{\partial^*E}\int_0^\e X(\eta,x)\cdot \nu_E(x)\,d\eta\,d\mathcal H^{N-1}(x)\\
&=\int_{\partial^*E}\lim_{\e\to0} X(\eta,x)\cdot \nu_E(x)\,d\eta\,d\mathcal H^{N-1}(x)\\
&=\int_{\partial^*E}\theta(\nu_E(x))\,d\mathcal H^{N-1}(x)
\end{aligned}
\]
and this ends the proof of Theorem \ref{main1}.\qed

\begin{remark}
We remark that if $E$ is a $C^{1,1}$-regular set in $\Omega$ then the computation of the pointwise limit is easier. Indeed, for such sets the following geometric property holds true (for details see \cite[Section I.2]{W}): there exists $r>0$ such that the map 
\[
\Psi_r \colon \partial E \times [0,r] \to \{y\in E^c : d(y,\partial E)\le r\}, \quad \Psi_r(x)=x+t\nu_E(x)
\]
is a $C^{1,1}$-diffeomorphism. Thus, performing change of variable $x=y-\e z$ we have 
\[
\begin{aligned}
\mathscr F_\e(E)&=\frac{1}{\e}\int_{\{y\in E^c : d(y,\partial E)\le \e\}}f\left(\frac{1}{\e^N}\int_EG\left(\frac{y-x}{\e}\right)dx\right)\,dy\\
&=\frac{1}{\e}\int_{\{y\in E^c : d(y,\partial E)\le \e\}}f\left(\int_{\frac{y-E}{\e}}G(z)dz\right)\,dy\\
&=\frac{1}{\e}\int_{\Psi_\e(\partial E \times [0,\e])}f\left(\int_{\frac{y-E}{\e}}G(z)dz\right)\,dy.
\end{aligned}
\]
For any $(x,t) \in \partial E \times [0,\e]$ let $J_\e(x,t)=|\det D\Psi_\e(x)|$. Then, using also $t=\e s$, 
\[
\begin{aligned}
\mathscr F_\e(E)&=\frac{1}{\e}\int_{\partial E}\int_0^\e f\left(\int_{\frac{x-E}{\e}+\frac{t}{\e}\nu_E(x)}G(z)dz\right)\,J_\e(x,t)\,d\mathcal H^{N-1}(x)\,dt\\
&=\int_{\partial E}\int_0^1 f\left(\int_{\frac{x-E}{\e}+s\nu_E(x)}G(z)dz\right)\,J_\e(x,\e s)\,d\mathcal H^{N-1}(x)\,ds.
\end{aligned}
\]
Since the regularity of $E$ we have 
\[
\lim_{\e\to0}J_\e(x,\e s)=1
\]
from which, applying again \eqref{conv}, 
\[
\lim_{\e\to0}\mathscr F_\e(E)=\int_{\partial E} \int_0^1f\left(\int_{\{v\cdot \nu_E(x)\ge t\})}G(z)\,dz\right)\,dt\,d\mathcal H^{N-1}(x)=\mathscr F(E).
\]
\end{remark}

\section{The lower estimate}

In this section we will prove our second main result, that is Theorem \ref{main2}. First of all at any $x\in D$ we let 
\[
\nu_h(x)=\frac{(-\nabla u_h(x),1)}{\sqrt{1+|\nabla u_h(x)|^2}}.
\]
It turns out that $\nu_h(x)$ is the exterior unit normal to $\partial^*E_h$ at $(x,u_h(x))$. For any $\eta>0$ small enough 
\[
D^\eta=\{x \in D : d(x,\partial D)>\eta\}.
\]
It turns out that $D^\eta \nearrow D$ in $L^1$ as $\eta\to 0^+$. If $z\in \R^N$ we will use the notation $z=(\bar z,z^N)$. We now divide the proof into several steps.
\\
\\
\textsc{Step 1:} We claim that for any $\sigma>0$, for any $x\in D^{3\sigma}$ and for any $h\in \N$ with $\e_h<\sigma$ we have 
\begin{equation}\label{step1}
\overline{B_2^{N-1}(0)}\subset \frac{x-D^\sigma}{\e_h}.
\end{equation}
Indeed, $x\in D^{3\sigma}$ means that $\overline{B^{N-1}_{2\sigma}(x)}\subset D^\sigma$. If now $z\in \R^{N-1}$ and $|z| \le 2$ then $|x-\e_hz-x| \le 2\e_h<2\sigma$ which implies that $x-\e_hz\in D^\sigma$ and then \eqref{step1}.
\\
\\
\textsc{Step 2:} For any $x\in D^{3\sigma}$, $s\in [0,1]$ and $\xi \in \R^{N-1}$ we let 
\[
a_h(x,s,\xi)=\frac{u_h(x)-u_h(x+\e_hs\overline{\nu_h(x)}-\e_h\xi)}{\e_h}+s\nu_h(x)^N.
\]
We claim that 
\begin{equation}\label{step2}
\lim_{h\to+\infty}a_h(x,s,\xi)=\nabla u(x) \cdot (\xi -s \overline{\nu(x)})+s\nu(x)^N.
\end{equation}
Indeed, 
\[
\begin{aligned}
&\frac{u_h(x+\e_hs\overline{\nu_h(x)}-\e_h\xi)-u_h(x)}{\e_h} \\
& = \frac{1}{\e_h} \int_0^{\e_h} \frac{d}{dt}u_h(x+t(s\overline{\nu_h(x)}-\xi)) \, dt \\
& = \frac{1}{\e_h} \int_0^{\e_h} \nabla u_h(x+t(s\overline{\nu_h(x)}-\xi)) \cdot (s\overline{\nu_h(x)}-\xi) \, dt \\
&= \frac{1}{\e_h} \int_0^{\e_h}  (\nabla u_h(x+t(s\overline{\nu_h(x)}-\xi))-\nabla u(x+t(s\overline{\nu_h(x)}-\xi))) \cdot (s\overline{\nu_h(x)}-\xi) \, dt\\
& \qquad \qquad +\frac{1}{\e_h} \int_0^{\e_h} (\nabla u(x+t(s\overline{\nu_h(x)}-\xi))-\nabla u(x+t(s\overline{\nu(x)}-\xi))) \cdot (s\overline{\nu_h(x)}-\xi) \, dt \\
&\qquad \qquad \quad +\frac{1}{\e_h} \int_0^{\e_h} \nabla u(x+t(s\overline{\nu(x)}-\xi)) \cdot (s\overline{\nu_h(x)}-\xi) \, dt=:I_1+I_2+I_3.
\end{aligned}
\]
Concerning the first integral, we have
\[
\begin{aligned}
|I_1|\le (s+|\xi|)\|\nabla u_h - \nabla u\|_{\infty} \rightarrow 0 \,\,\, \textnormal{as $h \rightarrow +\infty.$}
\end{aligned}
\]
On the other hand, if $L$ is the Lipschitz constant of $\nabla u$, we get
\[
I_2 \le L(s+|\xi|)\|\overline{\nu_h}-\overline{\nu}\|_\infty \rightarrow 0 \,\,\, \textnormal{as $h \rightarrow +\infty.$}
\]
Finally, for the third integral, let $g(t)= \nabla u(x+t(s\overline{\nu(x)}-\xi))$. Then $g$ is continuous, hence
\[
\lim_{h\to+\infty}\frac{1}{\e_h}\int_0^{\e_h} g(t) \, dt=g(0)
\]
from which 
\[
\lim_{h\to+\infty}I_3=\nabla u(x) \cdot (s \overline{\nu(x)}-\xi)
\]
as claimed.
\\
\\
\textsc{Step 3:}  Let $M=\sup_h\|u_h\|_\infty$ and let $\sigma \in (0,M/2)$. We claim that for any $h\in \N$ with $\e_h<\sigma$ it holds 
\begin{equation}\label{step3}
\mathscr F_{\e_h}(E_h)\ge \int_{D^{3\sigma}}\int_0^1f\left(\int_{\overline{B_1^{N-1}(0)}}\int_{a_h(x,s,\xi)}^1G(\xi,\eta) \, d \eta \,d \xi\right)ds\,\sqrt{1+|\nabla u_h(x)|^2}\,dx.
\end{equation}
Indeed, first of all notice that
\[
\{z \in E_h^c : B_{\e_h}(z) \cap E_h \ne \emptyset\} \supset \{(x-r\overline{\nu_h(x)},u_h(x)+r\nu_h(x)^N) : x\in D^\sigma,\,r \in (0,\e_h)\},
\]
As a consequence, 
\[
\mathscr F_{\e_h}(E_h) \ge \frac{1}{\e_h} \int_{D^{3\sigma}} \int_0^{\e_h} f \left (G_{\e_h}\ast \chi_{E_h}(x-r\overline{\nu_h(x)},u_h(x)+r\nu_h(x)^N)  \right ) \, dr \, \sqrt{1 + |\nabla u_h(x)|^2} \, dx.
\]
We concentrate now on the term $G_{\e_h}\ast \chi_{E_h}(x-r\overline{\nu_h(x)},u_h(x)+r\nu_h(x)^N)$ and we rewrite it in a suitable way by performing some changes of variables. First of all, by noticing that $E_h = \{(z,w) \in D \times \R : 0\le w \le u_h(z)\}$ we have
\[
\begin{aligned}
&G_{\e_h}\ast \chi_{E_h}(x-r\overline{\nu_h(x)},u_h(x)+r\nu_h(x)^N)\\
 & \ge \int_{D^\sigma} \frac{1}{\e_h^N} \int_0^{u_h(z)} G \left (\frac{x -r \overline{\nu_h(x)} - z}{\e_h}, \frac{u_h(x) + r\nu_h(x)^N - w}{\e_h} \right ) \, dw \, dz.
\end{aligned}
\]
We now perform the change of variables in the following order:
\[
\eta = \frac{u_h(x) + r\nu_h(x)^N  - w}{\e_h}, \quad \xi = \frac{x + r \overline{\nu_h(x)} - z}{\e_h}.
\]
We obtain 
\[
G_{\e_h}\ast \chi_{E_h}(x-r\overline{\nu_h(x)},u_h(x)+r\nu_h(x)^N)\ge\int_{\frac{x+r\overline{\nu_h(x)}-D^\sigma}{\e_h}}\int_{a_h(x,r/\e_h,\xi)}^{\frac{u_h(x)+r\nu_h(x)^N}{\e_h}}  G(\xi,\eta) \, d \eta \,d \xi.
\]
Recalling that $f$ is non-decreasing and operating the change of variable $r=\e_hs$ we arrive to 
\[
\mathscr F_{\e_h}(E_h) \ge \int_{D^{3\sigma}} \int_0^1 f \left( \int_{\frac{x-D^\sigma}{\e_h}+s\overline{\nu_h(x)}}\int_{a_h(x,s,\xi)}^{\frac{u_h(x)}{\e_h}+s\nu_h(x)^N}  G(\xi,\eta) \, d \eta \,d \xi\right) \, ds \, \sqrt{1 + |\nabla u_h(x)|^2} \, dx.
\]
Now, since \eqref{step1} we deduce that for any $x\in D^{3\sigma}$ and for any $s\in [0,1]$ 
\[
\frac{x-D^\sigma}{\e_h}+s\overline{\nu_h(x)} \supset \overline{B_2^{N-1}(0)}+s\overline{\nu_h(x)}\supset \overline{B_1^{N-1}(0)}.
\]
Moreover, using $\sigma<M/2$ we get also 
\[
\frac{u_h(x)}{\e_h}+s\nu_h(x)^N>1.
\]
As a consequence, recalling that $G$ is supported on $\overline{B_1^N(0)}$ we obtain \eqref{step3}.
\\
\\
\textsc{Step 4:} Passing to the limit as $h\to +\infty$ in \eqref{step3}, using Fatou's Lemma \eqref{step2} and the Dominated convergence Theorem we obtain 
\[
\begin{aligned}
\liminf_{h \rightarrow +\infty} &\mathscr{F}_{\e_h}(E_h)\\
& \ge \, \int_{D^{3\sigma}} \int_0^1 f \left ( \liminf_{h\to+\infty} \int_{\overline{B_1^{N-1}(0)}} \int_{a_h(x,s,\xi)}^1 G(\xi, \eta) \, d \eta \, d \xi\right ) ds \sqrt{1 +|\nabla u(x)|^2} \, dx\\
& = \int_{D^{3\sigma}} \int_0^1 f \left (  \int_{\overline{B_1^{N-1}(0)}} \int_{\nabla u(x) \cdot (\xi-s\overline{\nu(x)})+s\nu(x)^N}^1G(\xi, \eta) \, d \eta \, d \xi\right ) ds \,\, \sqrt{1 + |\nabla u(x)|^2} \, dx.
\end{aligned}
\]
By the arbitrariness of $\sigma$ small we get 
\[
\begin{aligned}
\liminf_{h \rightarrow +\infty}& \mathscr{F}_{\e_h}(E_h)\\
&\ge \int_D \int_0^1 f \left (  \int_{\overline{B_1^{N-1}(0)}} \int_{\nabla u(x) \cdot (\xi-s\overline{\nu(x)})+s\nu(x)^N}^1G(\xi, \eta) \, d \eta \, d \xi\right ) ds \,\, \sqrt{1 + |\nabla u(x)|^2} \, dx.
\end{aligned}
\]
\\
\textsc{step 5:} We conclude the proof showing that 
\[
\int_D \int_0^1 f \left (  \int_{\overline{B_1^{N-1}(0)}} \int_{\nabla u(x) \cdot (\xi-s\overline{\nu(x)})+s\nu(x)^N}^1G(\xi, \eta) \, d \eta \, d \xi\right ) ds \,\, \sqrt{1 + |\nabla u(x)|^2} \, dx = \mathscr{F}(E).
\]
First of all, we notice that 
\[
\eta=\nabla u(x) \cdot (\xi-s\overline{\nu(x)})+s\nu(x)^N=\nabla u(x) \cdot \xi+s\sqrt{1+|\nabla u(x)|^2}
\]
is the equation of an affine hyperplane in $\R^N$ orthogonal to $\nu(x)$ whose distance from the origin is 
\[
\frac{s\sqrt{1+|\nabla u(x)|^2}}{\sqrt{1+|\nabla u(x)|^2}}=s.
\]
As a consequence, 
\[
 \int_{\overline{B_1^{N-1}(0)}} \int_{\nabla u(x) \cdot (\xi-s\overline{\nu(x)})+s\nu(x)^N}^1G(\xi, \eta) \, d \eta \, d \xi=\int_{\{z\cdot\nu_E(x,u(x))\ge s\}} G(z) \, dz
\] 
from which 
\[
\begin{aligned}
\int_D \int_0^1& f \left (  \int_{\overline{B_1^{N-1}(0)}} \int_{\nabla u(x) \cdot (\xi-s\overline{\nu(x)})+s\nu(x)^N}^1G(\xi, \eta) \, d \eta \, d \xi\right ) ds \,\, \sqrt{1 + |\nabla u(x)|^2} \, dx \\
&=\int_D \int_0^1 f \left ( \int_{\{z\cdot\nu_E(x,u(x))\ge s\}} G(z) \, dz\right ) ds \,\, \sqrt{1 + |\nabla u(x)|^2} \, dx\\
&= \int_{\partial E} \int_0^1 f \left ( \int_{\{z\cdot\nu_E(y)\ge s\}} G(z) \, dz\right ) ds \, d\mathcal H^{N-1}(y)=\mathscr{F}(E)
\end{aligned}
\]
and the proof is complete.
\qed

\section{$L^1$-lower semicontinuity of $\mathscr F$}
 
We are going to prove Theorem \ref{main3}. It is well known (see for instance \cite[Thm.\,5.14]{AFP-libro}) that is sufficient to check that the positively one-homogeneous extension of $\theta$ given by 
\[
\tilde \theta(v)=\left\{\begin{array}{ll}
\displaystyle |v|\theta\left(\frac{v}{|v|}\right)\,dt & \text{if $v\ne 0$},
\\
0 & \text{if $v=0$},
\end{array}\right.
\]
is convex. First of all, by direct computation for each $v\in \R^N$ with $v\ne 0$ we have 
\[
\theta\left(\frac{v}{|v|}\right)= \int_0^1 f\left(\int_{\{z \cdot v\ge |v|t\}}G(z)\,dz\right)\,dt\stackrel{|v|t=s}{=}\frac{1}{|v|}\int_0^{|v|} f\left(\int_{\{z \cdot v\ge s\}}G(z)\,dz\right)\,ds
\]
from which we obtain
\[
\tilde \theta(v)=\left\{\begin{array}{ll}
\displaystyle \int_0^{|v|} f\left(\int_{\{z \cdot v\ge s\}}G(z)\,dz\right)\,ds & \text{if $v\ne 0$},
\\
0 & \text{if $v=0$}.
\end{array}\right.
\]
Now it is east to see that $\tilde\theta$ is convex. Indeed, since $f$ is convex there exist $(\alpha_h),(\beta_h)$ such that 
\[
\text{$f=\lim_{h\to+\infty}f_h$ uniformly on compact sets, where $f_h(t)=\alpha_ht+\beta_h$}.
\]
For any $h\in \N$ let 
\[
\tilde \theta_h(v)=\left\{\begin{array}{ll}
\displaystyle \int_0^{|v|} f_h\left(\int_{\{z \cdot v\ge s\}}G(z)\,dz\right)\,ds & \text{if $v\ne 0$},
\\
0 & \text{if $v=0$}.
\end{array}\right.
\]
Since $f_h\to f$ uniformly on $[0,1]$ we can say that $\tilde \theta_h \to \tilde \theta$ pointwise. In order to conclude it is sufficient to show that $\tilde\theta_h$ is convex. For any $v\ne 0$ we let $\hat v=\frac{v}{|v|}$. Then  
\[
\begin{aligned}
\tilde\theta_h (v)&=\alpha_h\int_0^{|v|} \int_{\{z \cdot v\ge s\}}G(z)\,dz\,ds+\beta_h|v|\\
&=\alpha_h\int_0^{|v|} \int_{\{\bar z \cdot v=0\}}\int_{s/|v|}^{+\infty}G(\bar z+t\hat v)\,d\bar z\,dt\,ds+\beta_h|v|\\
&=\alpha_h\int_0^{+\infty} \int_{\{\bar z \cdot v=0\}}\int_0^{t|v|}G(\bar z+t\hat v)\,ds\,dt\,d\bar z+\beta_h|v|\\
&=\alpha_h|v|\int_0^{+\infty} \int_{\{\bar z \cdot v=0\}}tG(\bar z+t\hat v)\,dt\,d\bar z+\beta_h|v|\\
&=\alpha_h\int_{\{\bar z \cdot v\ge 0\}}G(z)z\cdot v\,dz+\beta_h|v|\\
&=\frac{\alpha_h}{2}\int_{\R^N}G(z)|z\cdot v|\,dz+\beta_h|v|
\end{aligned}
\]
where the last equality follows since $G$ is even. Notice that the last expression is convex in $v$ and this ends the proof.

\section{Some examples}

In this section we characterize the limit functional $\mathscr F$ in some interesting cases. 

\subsection{$G$ radially symmetric}
Assume $G(z)=g(|z|)$ for some $g \colon [0,+\infty) \to \R$. Take $\nu\in \mathbb S^{N-1}$ and $t\ge 0$. Notice that the quantity 
\[
\int_{\{z\cdot \nu\ge t\}}G(z)\,dz
\]
does not depend on $\nu$. Take now $E\in \mathcal P_N$ and $x\in \partial^*E$. We have 
\[
\int_0^1f\left(\int_{\{z\cdot \nu_E(x)\ge t\}}G(z)\,dz\right)\,dt=c
\]
where $c$ is a constant that depends only on $N,f$ and $G$. Then 
\[
\mathscr F(E)=c\,\mathcal H^{N-1}(\partial^*E).
\]

\subsection{The case $f(t)=t$}
When $f$ is the identity function for any $E\in \mathcal P_N$ and for any $x\in \partial^*E$ we have
\[
\begin{aligned}
\theta(\nu_E(x))&=\int_0^1\int_{H_{\nu_E(x)}+t\nu_E(x)}G(z)\,dz\,dt\\
&=\int_0^1 \int_{\{\bar z\cdot \nu_E(x)=0\}}\int_t^1G(\bar z+s\nu_E(x))\,ds\,d\bar z\,dt\\
&=\int_0^1 \int_{\{\bar z\cdot \nu_E(x)=0\}}\int_0^sG(\bar z+s\nu_E(x))\,dt\,d\bar z\,ds\\
&=\int_0^1 \int_{\{\bar z\cdot \nu_E(x)=0\}}s\,G(\bar z+s\nu_E(x))\,d\bar z\,ds\\
&=\int_{H_{\nu_E(x)}}G(z)z\cdot \nu_E(x)\,dz\\
&=\frac{1}{2}\int_{\R^N}G(z)|z\cdot \nu_E(x)|\,dz.
\end{aligned}
\]
Then the limit $\mathscr F$ is given by 
\[
\mathscr F(E)=\frac{1}{2}\int_{\partial^*E}\int_{\R^N}G(z)|z\cdot \nu_E(x)|\,dz\,d\mathcal H^{N-1}(x).
\]
This is in accordance to \cite{P}.
\begin{remark}
If $N>1$ and $G$ is radially symmetric we have, if $g \colon [0,+\infty) \to \R$ is such that $G(z)=g(|z|)$, 
\[
\begin{aligned}
\frac{1}{2}\int_{\R^N}G(z)|z\cdot \nu_E(x)|\,dz&=\frac{1}{2}\int_{\R^N}g(|z|)|z\cdot \nu_E(x)|\,dz\\
&=\frac{1}{2}\int_0^{+\infty}\int_{\mathbb S^{N-1}}g(r)r|\xi \cdot \nu_E(x)|\,dr\,d\mathcal H^{N-1}(\xi)\\
&=|B^{N-1}_1(0)|\int_0^{+\infty}g(r)r\,dr\\
&=\frac{|B^{N-1}_1(0)|}{\mathcal H^{N-1}(\mathbb S^{N-1})}\int_{\R^N}G(z)|z|\,dz
\end{aligned}
\]
since it is well known that for any $\nu\in \mathbb S^{N-1}$ it holds
\[
\frac{1}{2}\int_{\mathbb S^{N-1}}|\xi \cdot \nu|\,d\mathcal H^{N-1}(\xi)=|B^{N-1}(0)|.
\]
We thus deduce that 
\[
\mathscr F(E)=c_{N,G}\,\mathcal H^{N-1}(E), \quad c_{N,G}=\frac{|B^{N-1}_1(0)|}{\mathcal H^{N-1}(\mathbb S^{N-1})}\int_{\R^N}G(z)|z|\,dz.
\]
This is in accordance to \cite{EO}.
\end{remark}

\section*{Acknowledgments}
We thank Irene Fonseca which proposed us this very nice problem.

\end{document}